Fernando Zalamea
(Universidad Nacional de Colombia)

# TWO NEW GESTURES ON PEIRCE'S CONTINUUM AND THE EXISTENTIAL GRAPHS

## 1. Introduction

In mathematics, the many blends of abstraction and concretion are ubiquitous. The connections can be seen as *gesture inversions* between opposites: universal/particular, negative/positive, one/multiple, topological/logical, etc. In fact, we have available many beautiful, concrete gestures where the transitions become technically definable: (1) Peirce's existential graphs, where the recto (positive) and the verso (negative) of the sheet of assertion allow the canonical transit of connectives (iteration/deiteration rule) and quantifiers (normal forms associated to the line of identity) (Zalamea 2012a), and where a multiplication of layers allows to connect modal logic (possible worlds, viewed as many sheets of assertion) with complex variable (the many sheets of a Riemann surface) (Zalamea 2010); (2) Poincaré's homotopy, where an understanding of paths in space opens-up the possibility to study equivalences, twists and inversions of paths, codified in homotopy groups (over a singular, base point) or homotopy groupoids (over plural, many points) (Zalamea 2012b); (3) Grothendieck's toposes, where, over a site (that is, a category with a family of covers well behaved = Grothendieck topology), one can define sheaves which connect local and global properties, and where a pulling-back of a canonical truth-value provides a universal algebra of multiple truth-values (subobject classifier), where the logic, always intuitionistic, ceases in general to be classical (Zalamea 2012b).

      In these examples, as in many others, we are looking at the very hand gesturing of the mathematician –traceable concrete residues– evolving in space. Some movements are easily described: *(A) pulling back-and-forth*: (1) graphs iterated and deiterated through

cuts, (2) paths twisted or inverted, (3) truth-values extended along the subobject classifier; *(B) plasticizing and deforming*: (1) an actual sheet projected/injected in many possible sheets, (2) a path deformed homotopically into another, (3) a sheaf contextually described by its synthetic gluing behavior following Yoneda's lemma; *(C) ascending and descending*: (1) a graph elevated through many layers of sheets of assertion to be further projected or actualized in a bottom one, (2) a topological space abstracted through the sequence of homotopy groups to be further characterized by that possible abstraction (Poincaré's Conjecture, in the case $S^3 \subseteq R^4$), (3) a site elevated in an arbitrary category to be further projected both in number (arithmetical toposes) and in space (geometrical toposes). Notice that gestures *(A)-(C)* also abound at the hands of an orchestra conductor, nothing really surprising, since mathematics and music have always been very close (Mazzola 2002).

In this article, we present two gestures corresponding to two profound new understandings of Peirce's Continuum (Vargas 2015) and Peirce's Existential Graphs (Oostra 2010). Vargas and Oostra have revolutionized Peirce's mathematical studies, thanks to a first complete model for Peirce's continuum provided by Vargas, and thanks to the emergence of intuitionistic existential graphs provided by Oostra. We will show below how these careful mathematical constructions can be encrypted in very simple gestures, along the lines *(A)-(C)* just mentioned. Some nice drawings by a young mathematician and artist, Angie Hugueth, reveal forcefully the mathematical gestures lying at the heart of Vargas and Oostra constructions.

## 2. Vargas' Complete Model for Peirce's Continuum

Peirce's continuum is a concept invented and explored by Peirce between 1890 and 1905 (Moore 2010), governed by some pretty strong requirements: (i) genericity and supermultitudeness, (ii) reflexivity and inextensibility, (iii) modality and plasticity (Zalamea 2012a). Some stringent consequences of these conditions are (i) that the size of the continuum cannot be determined (genericity) and thus exceeds any given infinite cardinal (supermultitudeness), (ii) that any non-void part is similar to the whole (reflexivity) and thus the continuum cannot contain points without further elements (inextensibility), and, finally, (iii) that all germs of possibility should be available in the con-tinuum (modality) and thus should include many logical layers in its construction (plasticity). For more than a century, an objectualization of Peirce's continuum concept, that is, a unifying structure that would include simultaneously (i)-(iii), was considered impossible by all Peirce scholars, either philosophers, logicians, or mathematicians. Nevertheless, for the great benefit of Peirce's studies, my student Francisco Vargas did not believe much in our impossibility claims, and, by a single stroke of genius has proved all us wrong (model discovered around 2012, published 2015).

As happens often in mathematics, the construction is as simple and profound as possible (following Grothendieck, simplicity should be mandatory when delving in mathematical depths). Just begin with iterations of the usual (cantorian) real line: for a non-null ordinal $\alpha$ consider the set $(C_\alpha, <_\alpha)$ of $\alpha$-real sequences with the lexicographical order. Then, extend the construction and define $(C, <)$ by taking the union of all $(C_\alpha)$'s over all ordinals and by taking $<$ as extension of the $(<_\alpha)$'s. Now, the core of the construction appears: for $x, y$ in $C$ define a new order relation E by *inversion*: $x$ E $y$ holds if $\text{dom}(y) < \text{dom}(x)$ and $x$ coincides with $y$ over $\text{dom}(y)$. In this way, the gigantic ordinal iterated tree of sequences is turned upside down. Every element in $C$ can then look at its descendants: for $x$ in $C$, we call such a collection $M_x = \{y : y \text{ E } x\}$ the monad associated to $x$. The main result of Vargas' construction (very simple proof!) is that, for any $x$, the monad is isomorphic to the whole ($M_x \approx C$). For details, see (Vargas 2015) and his other forthcoming articles. Genericity and supermultitudeness are obtained at once, due to the iteration along the whole proper class of ordinals. Reflexivity corresponds to the main result ($M_x \approx C$), and inextensibility follows at once from reflexivity. Modality is a little more hidden, but one can easily imagine how, below any given sequence, its descendants cover the realm of possibilities. Plasticity can be found at every level: genericity of size, isomorphism of the monads, possible descendants. In fact, Vargas' model reveals an outstanding richness, truly elastic in all technical, conceptual and metaphorical perspectives.

Iteration and inversion are the main clues of Vargas' model. It is clear that these ideas enter into the main gestural forces *(B)-(C)* stated at our *Introduction*. Vargas plasticizes and deforms, ascends and descends. In 2017, I tried to explain these simple ideas *with ges-* *tures only* to my very young student Angie Hugueth (just in second year mathematics): draw a horizontal line to represent the first stage (reals), then on any given point draw a vertical line (second stage), continue iterating indefinitely the lines, and finally, turn the structure upside down. Some days later, Angie came out with a nice drawing ('Gesto de Vargas'), where all is encompassed in a single image (see *Figure 1*). In red, the hands produce the basis and the iterations. In blue, the curves produce the idea of an inverted structure (curls of a maelstrom). In this way, a simple gesture is able to relate forcefully philosophical goals, definitions, techniques, proofs, and layers of imagination.

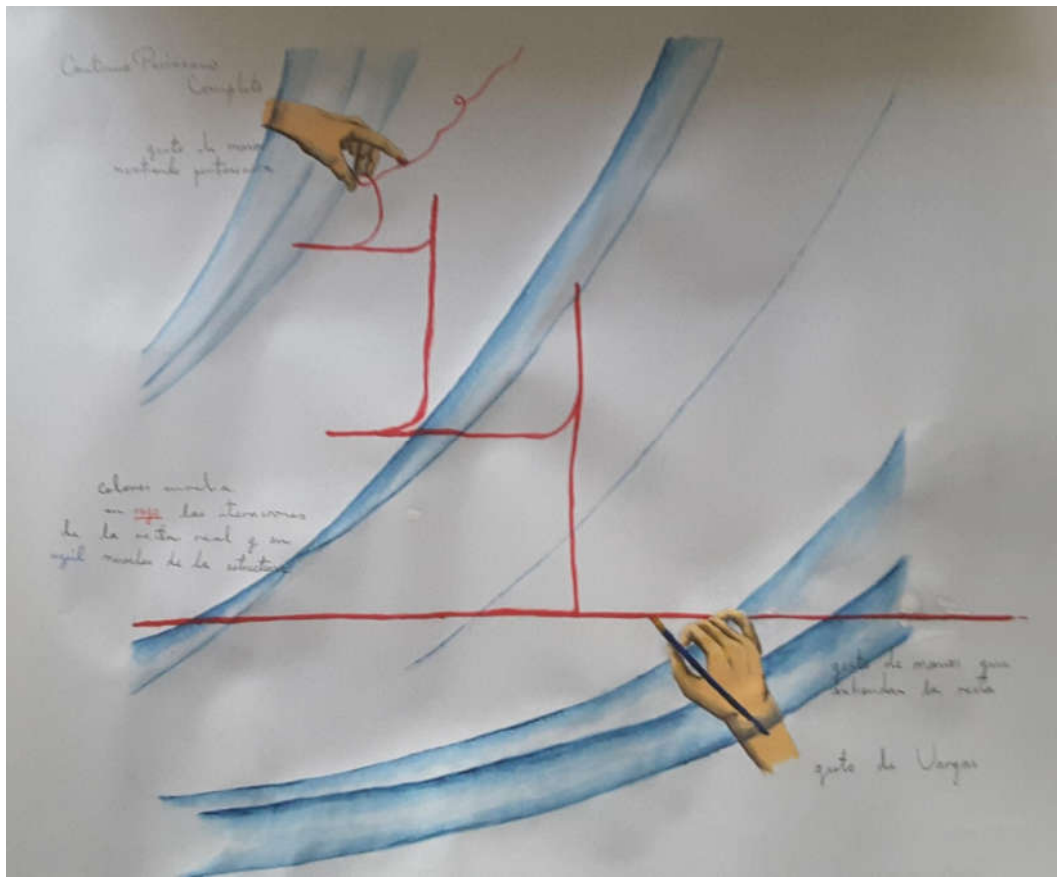

*Figure 1.* Angie Hugueth, *Gesto de Vargas* (2017)

### 3. Oostra's Intuitionistic Existential Graphs

Peirce invented the existential graphs also in the brilliant period 1890-1905, alongside with his explorations on the continuum. The existential graphs (Roberts 1963, Zeman 1964, Zalamea 2012a) provide an extremely original graphical logic of continuity, and constitute the only system of logic which, with the *same* rules, but varying the language, axiomatize simultaneously classical propositional calculus, first-order classical relational logic, and some modal calculi (including Lewis' systems). When you compare this situation with the extremely different axioms that, in Hilbert type systems, govern propositional and predicate calculi, one discovers that the existential graphs do reveal some *archetypes* of logical reasoning hidden in other presentations. The Alpha sheet of assertion, a continuous surface on which the graphs are marked, stands as an iconic reflection of real non-degenerate continuity (thirdness, following Peirce's cenopythagorean categories); Alpha characterizes the classical propositional calculus. The Beta line of identity, a continuous line which allows quantifying portions of reality, iconically reflects existence as degenerate continuity (secondness); Beta characterizes the classical first-order purely relational calculus. Finally, the Gamma half-cuts reflect openings to contingency and possibility (firstness); some forms of Gamma characterize Lewis' systems for modal propositional logic. For a short synthesis, see the *Table* below.

```
1. Signs.

Sheet of assertion:    blank generic sheet.              Icon:    □

Cuts:                  generic ovals detaching regions
                       in the sheet of assertion.        Icons:   ◯  ⬭
                                                                 (Alpha) (Gamma)

Line of identity:      generic line weaving relations
                       in the sheet of assertion.        Icon:    ―――
                                                                  (Beta)

Logical terms:         propositional and relational signs
                       marking the sheet of assertion.   Icons:   p, q, ...  R, S, ....

2. Illative Transformations of Signs.

Detaching Properties ("information zones").
Cuts can be nested but cannot intersect.
Identity lines can intersect other identity lines and all kinds of cuts.

Double cuts Alpha can be introduced or eliminated around any graph, whenever in
the "donut" region (gray) no graphs different from identity lines appear.            ◎

Transferring Properties ("information transmission").
Inside regions nested in an even number of cuts (Alpha or Gamma), graphs may be erased.
Inside regions nested in an odd number of cuts, graphs may be inserted.
Towards regions nested in a bigger number of cuts, graphs may be iterated.
Towards regions nested in a lower number of cuts, graphs may be deiterated.

3. Interpretation of Signs and Illative Transformations.

Blank sheet:           truth
Alpha cut:             negation
Juxtaposition:         conjunction
Line of identity:      existential quantifier
Gamma cut:             contingency (possibility of negation)

Double cut:            classical rule of negation  ($\neg\neg p \leftrightarrow p$)
Erasure and insertion: minimal rule of conjunction ($p \wedge q \rightarrow p$  and  $\neg p \rightarrow \neg(p \wedge q)$)
Iteration and deiteration: intuitionistic rule of negation as generic connective ($p \wedge \neg q \leftrightarrow p \wedge \neg(p \wedge q)$)
```

*Table*
*A quick view of the Existential Graphs systems:* ALPHA, BETA, GAMMA

I have had the chance of having fantastic students, which, then, with time, have become my teachers. Arnold Oostra is possibly the greatest student I have encountered, and, undoubtedly, he has now become the foremost international scholar on Peirce's mathematics. Around 2007, I was discussing with Oostra the question of why, since existential graphs are naturally connected to a logic of continuity, we had Alpha graphs corresponding to classical propositional calculus, instead of some *intuitionistic* existential graphs, more akin to continuity[1]. After many trials, Oostra arrived to a very simple and profound idea, and showed that, in fact, intuitionistic existential graphs could be constructed as another form of the archetypical Peircean methodology. Oostra kept again the *same* rules, but introduced new graphical signs for intuitionistic implication and disjunction (Oostra 2010), and proved, through extremely simple graphical means, the known correlations between classical and intuitionistic calculi. The core of the construction consists in comparing the sign of (i) an Alpha classical implication (through a discrete double cut) and (ii) an Alpha intuitionistic implication (through a continuous double cut: two cuts glued in one point). The rules of the existential graphs show then how to pass from (ii) to (i) (an intuitionistic proof is a classical one), but *not* viceversa. In other words, the strong logical differences between classicism and intuitionism[2] are captured in the existential graphs by processes of glueing and separation of nested cuts. Now, these ideas enter into the main gestural forces *(A)-(B)* stated at our *Introduction*. Oostra pulls back-and-forth, plasticizes and deforms. In this way, Oostra shows that the natural graphs akin to a logic of continuity, or topological logic, are in fact the intuitionistic existential graphs (Oostra 2010, 2011).

Again, in my 2017 beautiful session with Angie Hugueth, I tried to explain to her, by *hands only*, the achievements obtained by Oostra. It was very easy: just draw, with a single stroke, two nested circles, touching each other in one point, and then pull them apart; by logical constraints related to the Alpha rules, sometimes a back-and-forth of the hands may be obtained, sometimes not. With the imaginative freedom of youth, and her ability to entangle art and mathematics, Angie came up with a 'Gesto de Oostra' (*Figure 2*), which sums up in a single image all the many subtleties of classicism versus intuitionism. Here, the red curves describe the positive parts of the sheet (pair regions), the blue curves hint at the negative parts (odd regions), and the white curves capture the processes of glueing and separating the nested cuts.

---

[1] Intuitionism and topology were naturally close since Brouwer, and became more so with Tarski's proof that a complete semantics for intuitionism is given by the class of topological spaces, and with Lawvere's proof that the intrinsic logic of any elementary topos is intuitionistic.

[2] Law of excluded middle, proof by contradiction, De Morgan rules, provable in classical logic, not in intuitionistic logic; discrete semantics (Boolean algebras) for classical logic, continuous semantics (Heyting algebras, topological spaces) for intuitionistic logic, etc.

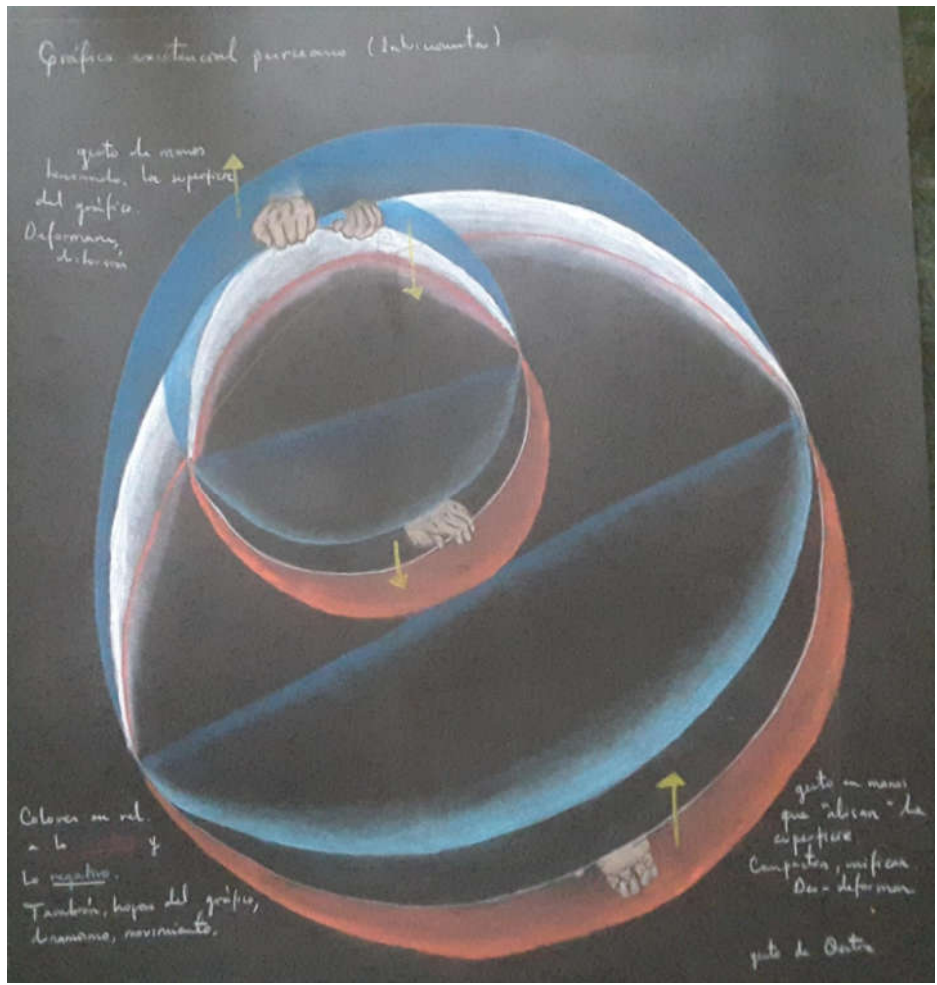

*Figure 2.* Angie Hugueth, *Gesto de Oostra* (2017)

## 4. Conclusions
I have no doubt that Francisco Vargas and Arnold Oostra achievements constitute, after Zeman's Thesis, the two single most important contributions to an understanding of Peirce's logic and mathematics provided in one hundred years. Vargas' model ends, once and for all, with the many fantastic speculations and controversies around Peirce's continuum. The model is simple, exact, convincing well beyond our doubts and prejudices. On another hand, Oostra's intuitionistic existential graphs show that a topological logic is best provided by Peirce's graphical methodology, and explain that some movements in the plane (related to Alpha rules: insertion/erasure, iteration/deiteration) are governed naturally by an underlying intuitionistic logic. Intuitionism, topology and visualization get thus strongly entangled. It is striking that such powerful achievements are so simple at the end, and that they can be captured by soft and precise gestures. This shows again the profound harmony between logic, mathematics and art, encrypted in some basic movements of the very human body.

## Acknowledgements
The main bulk of these ideas was presented at the 10th Anniversary of the *Associazione Pragma*, Università degli Studi di Milano, October 18 2017. I thank Rossella Fabbrichesi, Giovanni Maddalena and Rosa Calcaterra for the invitation, and Andrea Parravicini, Maria Regina Brioschi and Guido Baggio for the material organization. I had the opportunity to talk with my dear Italian colleagues several times in Milano (and Roma, a week before) around Peirce, mathematics and philosophy, always profiting from the great originality that the Italian school conveys on Peirce. My deep gratitude to all of them.